\documentclass[11pt]{amsart} 

\usepackage{amscd}
\usepackage{amssymb} 
\usepackage[text={6in,9in},centering]{geometry}

\newcommand{\ve}{\varepsilon}
\newcommand{\ddbar}{\sqrt{-1}\partial \overline{\partial}}
\newcommand{\nablabar}{\overline{\nabla}}
\newcommand{\pbar}{{\overline{p}}}
\newcommand{\qbar}{{\overline{q}}}
\newcommand{\rbar}{{\overline{r}}}
\newcommand{\ibar}{{\overline{i}}}
\newcommand{\jbar}{{\overline{j}}}
\newcommand{\kbar}{{\overline{k}}}
\newcommand{\lbar}{{\overline{l}}}

\newcommand{\sbar}{{\overline{s}}}

\newcommand{\abar}{{\overline{a}}}
\newcommand{\bbar}{{\overline{b}}}
\newcommand{\cbar}{{\overline{c}}}
\newcommand{\dbar}{{\overline{d}}}
\newcommand{\ov}[1]{\overline{#1}}

\newcommand{\Rm}{{\mathrm{Rm}}}
\newcommand{\ddt}{\frac{\partial}{\partial t}}
\newcommand{\tr}{{\mathrm{tr}_{\hat{g}}g}}

\begin{document}

\newcommand{\thmcounter}{section}
\newtheorem{theorem}{Theorem}[section]
\newtheorem{proposition}[theorem]{Proposition}
\newtheorem{claim}[theorem]{Claim}
\newtheorem{lemma}[theorem]{Lemma}
\newtheorem{corollary}[theorem]{Corollary}
\newtheorem{algorithm}[theorem]{Algorithm}
\newtheorem{question}[theorem]{Question}
\newtheorem{remark}[theorem]{Remark}
\numberwithin{equation}{section}

\title[Local estimates for the Chern-Ricci flow]{Local Calabi and curvature estimates for the Chern-Ricci flow$^{\dagger}$}

\author[M. Sherman]{Morgan Sherman}
\address{Department of Mathematics, California Polytechnic State University, San Luis Obispo, CA 93407}
\author[B. Weinkove]{Ben Weinkove}
\address{Department of Mathematics, Northwestern University, 2033 Sheridan Road, Evanston, IL 60208}
\thanks{$^{\dagger}$Supported in part by NSF grant DMS-1105373.  Part of this work was carried out while the second-named author was a member of the  mathematics department of the University of California, San Diego.}
\begin{abstract}
Assuming local uniform bounds on the metric for a solution of the Chern-Ricci flow, we establish local Calabi and curvature estimates using the maximum principle.   \end{abstract}

\maketitle

\section{Introduction}

Let $(M, \hat{g})$ be a Hermitian manifold.  The \emph{Chern-Ricci flow} starting at $\hat{g}$ is a smooth flow of Hermitian metrics $g=g(t)$ given by
\begin{equation}\label{crf0}
\ddt g_{i\ov{j}} = - R^C_{i\ov{j}}, \qquad g_{i\ov{j}}|_{t=0} = \hat{g}_{i\ov{j}},
\end{equation}
where $R^C_{i\ov{j}} := - \partial_i \partial_{\ov{j}} \log \det g$ is the \emph{Chern-Ricci} curvature of $g$.  If $\hat{g}$ is K\"ahler, then the Chern-Ricci flow coincides with the K\"ahler-Ricci flow.  

The Chern-Ricci flow was introduced by Gill \cite{G} and further investigated by Tosatti and the second-named author \cite{TW1, TW2}.  This flow has many of same properties as the K\"ahler-Ricci flow.  For example: on manifolds with vanishing first Bott-Chern class the Chern-Ricci flow converges to a Chern-Ricci flat metric  \cite{G};  on manifolds with negative first Chern class, the Chern-Ricci flow takes any Hermitian metric to the K\"ahler-Einstein metric \cite{TW1}; 
when $M$ is a compact complex surface and $\hat{g}$ is $\partial \ov{\partial}$-closed, the Chern-Ricci flow exists until either the volume of the manifold goes to zero or the volume of a curve of negative self-intersection goes to zero \cite{TW1}; if in addition $M$ is  non-minimal with nonnegative Kodaira dimension, the Chern-Ricci flow shrinks exceptional curves in finite time \cite{TW2} in the sense of Gromov-Hausdorff.  These results are closely analogous to results for the K\"ahler-Ricci flow \cite{Cao, FIK, TZ, SW0, SW}.

In this note, we establish local derivative estimates for solutions of the Chern-Ricci flow assuming local uniform bounds on the metric, generalizing our previous work \cite{ShW} on the K\"ahler-Ricci flow.   Our estimates are local, so we  work in a small open subset of  $\mathbb{C}^n$.  Write $B_r$ for the ball of radius $r$ centered at the origin in $\mathbb{C}^n$, and fix $T<\infty$.   We have the following result (see Section \ref{sectionprelim} for more details about the notation).

\begin{theorem} \label{maintheorem} Fix $r$ with $0< r<1$.  Let $g(t)$ solve the Chern-Ricci flow (\ref{crf0}) in a neighborhood of $B_r$ for $t\in [0,T]$.
Assume $N>1$ satisfies
\begin{equation} \label{uniformestimatemetric}
\frac{1}{N} \hat{g} \le g(t) \le N \hat{g} \qquad \textrm{on } B_r \times [0,T].
\end{equation}
Then there exist positive constants $C, \alpha, \beta$ depending only on $\hat{g}$ such that
\begin{enumerate}
\item[(i)] $\displaystyle{| \hat{\nabla} g|^2_{g} \le \frac{C N^{\alpha}}{r^2}}$ on $B_{r/2} \times [0,T],$ where $\hat{\nabla}$ is the Chern connection of $\hat{g}$.
\item[(ii)] $\displaystyle{|\emph{Rm}|_g^2 \le \frac{C N^{\beta}}{r^4}}$ on $B_{r/4} \times [0,T]$, for $\emph{Rm}$ the Chern curvature tensor of $g$.
\end{enumerate}
\end{theorem}

Note that the estimates are independent of the time $T$ and so the results holds also for time intervals $[0,T)$ or $[0,\infty)$.  The dependence of the constants on $\hat{g}$ is as follows: up to three derivatives of torsion of $\hat{g}$ and one derivative of the Chern curvature of $\hat{g}$ (see Remarks \ref{remark1} and \ref{remark2}).  We call the bound (i)  a \emph{local Calabi estimate} \cite{Ca} (see \cite{ZZ} for a similar estimate in the elliptic case).
 
 As a consequence of Theorem \ref{maintheorem}, we have local derivative estimates for $g$ to all orders:

\begin{corollary} \label{corollary}
With the assumptions of Theorem \ref{maintheorem}, for any $\ve>0$ with $0< \ve<T$, there exist constants $C_m$, $\alpha_m$ and $\gamma_m$ for $m= 1, 2,3,  \ldots$ depending only on $\hat{g}$ and $\ve$ such that 
$$| \hat{\nabla}^m_{\mathbb{R}} g|^2_{\hat{g}} \le \frac{C_m N^{\alpha_m}}{r^{\gamma_m}} \qquad \textrm{on } B_{r/8} \times [\ve,T],$$
where $\hat{\nabla}_{\mathbb{R}}$ is the Levi-Civita covariant derivative associated to $\hat{g}$.
\end{corollary}

Note that our assumption (\ref{uniformestimatemetric}) often holds for the Chern-Ricci flow on compact subsets away from a subvariety.  For example, this always occurs for the Chern-Ricci flow on a non-minimal complex surface of nonnegative Kodaira dimension \cite{TW1, TW2}.  It has already been shown by Gill \cite{G} that local derivative estimates exist using the method of Evans-Krylov \cite{E,K} adapted to this setting.  The purpose of this note is to give a direct maximum principle proof of Gill's estimates, and in the process identify evolution equations for the Calabi quantity $|\hat{\nabla} g|^2_g$ and the Chern curvature tensor $R_{i\ov{j}k\ov{l}}$, which were previously unknown for this flow.  In addition, we more precisely determine the form of dependence on the constants $N$ and $r$.  We anticipate that this may be useful, for example in generalizations of arguments of \cite{SW}.  

In the case when $\hat{g}$ is K\"ahler, so that $g(t)$ solves the K\"ahler-Ricci flow, the above result follows from results of the authors in \cite{ShW}. The more general case we deal with here leads to many more difficulties, arising from the torsion tensors of $g$ and $\hat{g}$.  For these reasons, our conclusions here are slightly weaker:  for example, we cannot obtain the small values ($\alpha=3$ and $\beta=8$)  in the estimates of (i) and (ii) that we achieved in \cite{ShW}.  

The second-named author thanks Valentino Tosatti and Xiaokui Yang for some helpful discussions.

\section{Preliminaries} \label{sectionprelim}

In this section we introduce the basic notions that we will be using throughout the paper.  We largely follow notation given in \cite{TW1}.  Given a Hermitian metric $g$ we write  $\nabla$ for the \emph{Chern connection} associated to $g$, which is characterized as follows.  Define Christoffel symbols $\Gamma_{ik}^l = g^{\sbar l} \partial_i g_{k\sbar}$.  Let $X = X^l \frac{\partial}{\partial z^l}$ be a vector field and let $a = a_k \, dz^k$ be a $(1,0)$ form.  Then 
\begin{equation}
	\nabla_i X^l = \partial_i X^l + \Gamma_{i r}^l X^r , 
	\quad \nabla_i a_j = \partial_i a_j - \Gamma_{i j}^r a_r .
\label{Chern connection}
\end{equation}
We can, in a natural way, extend $\nabla$ to act on any tensor.  Note that $\nabla$ makes $g$ parallel:  i.e. $\nabla g = 0$.  Similarly we let $\hat\nabla$ denote the Chern connection associated to $\hat g$.

Define the torsion tensor $T$ of $g$ by
\begin{equation}
	T_{ij}{}^k = \Gamma_{ij}^k - \Gamma_{ji}^k
\label{torsion}
\end{equation}
We note that $g$  is K\"ahler precisely when $T = 0$.  We write 
$T_{\ibar \jbar}{}^\kbar := \Gamma_{\ibar\jbar}^\kbar - \Gamma_{\jbar\ibar}^\kbar 
:= \overline{\Gamma_{ij}^k} - \ov{\Gamma_{ji}^k}$ for the components of the tensor $\overline{T}$.  
We lower and raise indices using the metric $g$.  For example, $T^{ij}{}_{k} = g^{\ov{a}i} g^{\ov{b}j} g_{k\ov{l}} T_{\ov{a} \ov{b}}{}^{\ov{l}}$.

We define the \emph{Chern curvature tensor} of $g$ to be the tensor written locally as  
\begin{equation}
	R_{i \jbar k}{}^l = - \partial_{\jbar} \Gamma_{ik}^l.
\label{curvature}
\end{equation}
Then
\begin{equation}
	R_{i \jbar k \lbar} = - \partial_i \partial_\jbar g_{k\lbar} 
			+ g^{\sbar r} \partial_i g_{k \sbar} \partial_\jbar g_{r \lbar} .
\label{curvature in terms of metric}
\end{equation}
where again we have lowered an index using the metric $g$.
Note that $\overline{R_{i \jbar k \lbar}} = R_{j \ibar l \kbar}$ holds.  

The commutation formulas for the Chern connection are given by
\begin{gather}
	[\nabla_i, \nabla_\jbar] X^l = R_{i \jbar k}{}^l X^k, 
		\quad [\nabla_i, \nabla_\jbar] \overline{X^k} 
		= - R_{i \jbar}{}^{\kbar}{}_{\lbar} \overline{X^l} \nonumber\\
	[\nabla_i, \nabla_\jbar] a_k = - R_{i \jbar k}{}^l a_l,  
		\quad [\nabla_i, \nabla_\jbar] \overline{a_l} 
		= R_{i \jbar}{}^{\kbar}{}_{\lbar} \overline{a_k}.
\label{curvature via difference in connections}
\end{gather}

Because $g$ is not assumed to be a K\"ahler metric the \emph{Bianchi identities} will not necessarily hold for $R_{i \jbar k \lbar}$.  However their failure to hold can be measured with the torsion tensor $T$ defined above:
\begin{align}
	& R_{i \jbar k \lbar} - R_{k \jbar i \lbar} = - \nabla_\jbar T_{i k \lbar} 
		\nonumber\\
	& R_{i \jbar k \lbar} - R_{i \lbar k \jbar} = - \nabla_i T_{\jbar  \lbar k} 
		\nonumber\\
	& R_{i \jbar k \lbar} - R_{k \lbar i \jbar} 
		= - \nabla_\jbar T_{i k \lbar} - \nabla_k T_{\jbar  \lbar i} 
		= - \nabla_i T_{\jbar  \lbar k} - \nabla_\lbar T_{i k \jbar} 
		\nonumber\\
	& \nabla_p R_{i \jbar k \lbar} - \nabla_i R_{p \jbar k \lbar} 
		= - T_{p i}{}^r R_{r \jbar k \lbar} 
		\nonumber\\
	& \nabla_\qbar R_{i \jbar k \lbar} - \nabla_\jbar R_{i \qbar k \lbar}
		= - T_{\qbar \jbar}{}^\sbar R_{i \sbar k \lbar}.
\label{Bianchi}
\end{align}
These identities are well-known (see \cite{TWY} for example).  Indeed,
it is routine to verify the first line, and the second and third lines follow directly from it.  Furthermore the fifth line follows directly from the fourth.  For the fourth line we calculate:
\[
	\nabla_p R_{i \jbar k}{}^l 
	= - \nabla_p (\partial_\jbar \Gamma_{ik}^l) 
	= - \partial_p \partial_\jbar \Gamma_{ik}^l 
		- \Gamma_{pr}^l \partial_\jbar \Gamma_{ik}^r 
		+ \Gamma_{pi}^r \partial_\jbar \Gamma_{rk}^l
		+ \Gamma_{pk}^r \partial_\jbar \Gamma_{ir}^l.
\]
Swapping the $p$ and $i$ indices, subtracting, and combining terms, we find
\[
	\nabla_p R_{i \jbar k}{}^l - \nabla_i R_{p \jbar k}{}^l 
	= - T_{pi}{}^r R_{r \jbar k}{}^l 
		+ \partial_\jbar \left( 
		\partial_i \Gamma_{pk}^l - \partial_p \Gamma_{ik}^l 
		+ \Gamma_{ir}^l \Gamma_{pk}^r 
		- \Gamma_{pr}^l \Gamma_{ik}^r
		\right).
\]
Now one checks that the quantity in parentheses vanishes.  

We define the Chern-Ricci curvature tensor $R^C_{i \jbar}$ by 
\begin{equation}
	R^C_{i \jbar} = g^{\lbar k} R_{i \jbar k \lbar} = - \partial_i \partial_{\jbar} \log \det g.
\label{Chern-Ricci curvature}
\end{equation}
Note that $\sqrt{-1} R^C_{i\ov{j}}dz^i \wedge dz^{\ov{j}}$ is a real closed (1,1) form.
We will suppose that $g = g(t)$ satisfies the \emph{Chern-Ricci flow}:
\begin{equation}
	\ddt g_{i \jbar} = - R^C_{i \jbar}, \quad g_{i\ov{j}}|_{t=0} = {\hat g}_{i\ov{j}}, 
\label{Chern-Ricci flow}
\end{equation}
for $t \in [0,T]$ for some fixed positive time $T$.  We will use $\hat{\nabla}$, $\hat{\Gamma}^l_{ik}$, $\hat{T}_{ik}{}^l$, $\hat{R}_{i\ov{j}k\ov{l}}$ etc to denote the corresponding quantities with respect to the metric $\hat{g}$.
Define a real (1,1) form $\omega=\omega(t)$ by $\omega = \frac{\sqrt{-1}}2 g_{i \jbar} dz^i \wedge dz^\jbar$ and similarly for $\hat\omega$.  From (\ref{Chern-Ricci flow}) we have that

\begin{equation}
	\omega = \hat\omega + \eta(t)
\label{omega = omegahat + eta}
\end{equation}
for a closed $(1,1)$ form $\eta$.  Hence
\begin{equation}
	T_{i k \lbar} = \hat T_{i k \lbar} .
\label{T = T hat}
\end{equation}
Here we raise and lower indices of $\hat{T}$ using the metric $\hat{g}$, in the same manner as for $g$ above.  Note that 
$T_{i k \lbar} = g_{r \lbar} T_{i k}{}^r = \partial_i g_{k \ov{l}} - \partial_k g_{i \ov{l}}$ and 
$\hat T_{i k \lbar} = \hat g_{r \lbar} \hat T_{i k}{}^r = \partial_i \hat{g}_{k \ov{l}} - \partial_k \hat{g}_{i \ov{l}}$.  

It is convenient to introduce the tensor $\Psi_{ik}{}^l = \Gamma_{ik}^l - \hat\Gamma_{ik}^l$.  We raise and lower indices of $\Psi$ using the metric $g$, and  write $\Psi_{\ov{i} \ov{k}}{}^{\ov{l}}$ for the components of $\ov{\Psi}$.
  We note here that $\Psi$ can be used to switch between the connections $\nabla$ and $\hat\nabla$.  For example given a tensor of the form $X_i{}^j$ we have 
\begin{equation}
	\nabla_p X_i{}^j - \hat\nabla_p X_i{}^j = - \Psi_{p i}{}^r X_r{}^j + \Psi_{p r}{}^j X_i{}^r . 
\label{nabla - nablahat = Psi}
\end{equation} 
Observe that
\begin{equation}
	\nabla_\jbar \Psi_{i k}{}^l = - R_{i \jbar k}{}^l + \hat R_{i \jbar k}{}^l.
\label{nablabar Psi}
\end{equation}
We write  $\Delta$ for the ``rough Laplacian'' of $g$, $\Delta = \nabla^{\qbar} \nabla_{\qbar}$, where $\nabla^\qbar = g^{\qbar p} \nabla_p$.  Finally note that we will write all norms $|\cdot|$ with respect to the metric $g$.

\section{Local Calabi estimate}

In this section we prove part (i) of Theorem \ref{maintheorem}.
We consider the Calabi-type \cite{Ca, Y} quantity  
\begin{equation}
	S := | \Psi |^2 = | \hat \nabla g |^2.
\label{definition of S}
\end{equation}
 
Our goal in this section is to uniformly bound $S$ on the set $B_{r/2}$, which we will do using a  maximum principle argument.  First we compute its evolution.  Calculate
\begin{align*}
	\Delta S ={} & g^{\qbar p} \nabla_p \nabla_{\qbar} 
		\left(g^{\abar i} g^{\bbar j} g_{k\cbar} \Psi_{ij}{}^k \overline{\Psi_{ab}{}^c} 
		\right) \\
	={} & g^{\qbar p} g^{\abar i} g^{\bbar j} g_{k \cbar} \nabla_p \left(
		\nabla_\qbar \Psi_{ij}{}^k \overline{\Psi_{ab}{}^c} + \Psi_{ij}{}^k 
		\overline{\nabla_q \Psi_{ab}{}^c } \right) \\
	={} & |\nablabar \Psi|^2 + |\nabla \Psi|^2 
		+ g^{\abar i} g^{\bbar j} g_{k \cbar} \Bigl( 
			\Delta\Psi_{ij}{}^k \overline{\Psi_{ab}{}^c} \\
		&+ \Psi_{ij}{}^k \overline{\left( 
			\Delta\Psi_{ab}{}^c 
			+ g^{\qbar p} R_{p \qbar a}{}^r \Psi_{rb}{}^c
			+ g^{\qbar p} R_{p \qbar b}{}^r \Psi_{ar}{}^c 
			- g^{\qbar p} R_{p \qbar r}{}^c \Psi_{ab}{}^r
			\right)} 
			\Bigr) \\
	={} & |\nablabar \Psi|^2 + |\nabla \Psi|^2 
		+ 2 \mathrm{Re} \left( 	( \Delta \Psi_{ij}{}^k ) \Psi^{ij}{}_k  \right) \\
		&+ (R_p{}^p{}_r{}^i \Psi^{rj}{}_k  
			+ R_p{}^p{}_r{}^j \Psi^{ir}{}_k  
			- R_p{}^p{}_k{}^r \Psi^{ij}{}_r  ) \Psi_{ij}{}^k.
\end{align*} 
From (\ref{nablabar Psi}) we have
\begin{equation}
	\Delta \Psi_{ij}{}^k = - \nabla^{\qbar} R_{i \qbar j}{}^k 
		+ \nabla^{\qbar} \hat{R}_{i \qbar j}{}^k.
\label{laplacian Psi}
\end{equation}

For the time derivative of $S$, first compute (cf. \cite{PSS} in the K\"ahler case), 
\begin{equation}
	\ddt \Psi_{ij}{}^k = \ddt \Gamma_{ij}^k = - \nabla_i (R^C)_j{}^k .
\label{evolution of Psi}
\end{equation}
Then
\begin{align*}
	\ddt S ={} & \ddt 
		\left( g^{\abar i} g^{\bbar j} g_{k\cbar} \Psi_{ij}{}^k \overline{\Psi_{ab}{}^c} 
		\right) \\
	={} & \left(\ddt g^{\abar i}\right) \Psi_{ij}{}^k \Psi_{\abar}{}^{j}{}_k 
		+ \left(\ddt g^{\bbar j}\right) \Psi_{ij}{}^k \Psi^i{}_{\bbar k} 
		+ \left(\ddt g_{k\cbar}\right) \Psi_{ij}{}^k \Psi^{i j \cbar} 
		+ 2\mathrm{Re} \left( \left(\ddt \Psi_{ij}{}^k\right) \Psi^{ij}{}_k \right) \\
	={} & (R^C)^{\abar i} \Psi_{ij}{}^k \Psi_{\abar}{}^{j}{}_k 
		+ (R^C)^{\bbar j} \Psi_{ij}{}^k \Psi^i{}_{\bbar k} 
		- (R^C)_{k\cbar} \Psi_{ij}{}^k \Psi^{i j \cbar} - 2\mathrm{Re} \left( (\nabla_i (R^C)_j{}^k) \Psi^{ij}{}_k \right).
\end{align*}
Therefore 
\begin{align*}
\left( \ddt - \Delta \right) S 
	={}& - |\nablabar \Psi|^2 - |\nabla \Psi|^2 
	+ \left( R^{\rbar i}{}_p{}^p - R_p{}^{p \rbar i} \right)
		\Psi_{ij}{}^k \Psi_{\rbar}{}^j{}_k 
	+ \left( R^{\rbar j}{}_p{}^p - R_p{}^{p \rbar j} \right)
		\Psi_{ij}{}^k \Psi^{i}{}_{\rbar k} \\
	&- \left( R_{k \rbar}{}_p{}^p - R_p{}^p{}_{k \rbar } \right)
		\Psi_{ij}{}^k \Psi^{i j \rbar} 
	- 2 \mathrm{Re} \left[ 
		\left( \nabla_i R_j{}^k{}_p{}^p + \Delta \Psi_{ij}{}^k \right) 
		\Psi^{i j}{}_k  \right].
\end{align*}
By (\ref{Bianchi}) we can re-write the terms involving a difference in curvature using the torsion tensor $T$.  For the term in square brackets we compute, using (\ref{laplacian Psi}) and again (\ref{Bianchi}) that 
\begin{align*}
	\nabla_i R_j{}^k{}_p{}^p + \Delta \Psi_{ij}{}^k 
		=& \nabla_i \left( 
			R_{p}{}^p{}_j{}^k + \nabla_j T^{p k}{}_p + \nabla^p T_{pj}{}^k
			\right) 
			- \nabla^\qbar R_{i \qbar j}{}^k 
			+ \nabla^\qbar \hat{R}_{i \qbar j}{}^k \\
		=& \left( 
			\nabla_p R_{i}{}^p{}_j{}^k - T_{ip}{}^r R_{r}{}^p{}_j{}^k 
			+ \nabla_i \nabla_j T^{p k}{}_p + \nabla_i \nabla^p T_{pj}{}^k
			\right) 
			- \nabla^\qbar R_{i \qbar j}{}^k + \nabla^\qbar \hat{R}_{i \pbar j}{}^k \\
		=& - T_{ip}{}^r R_r{}^p{}_j{}^k 
			+ \nabla_i \nabla_j T^{p k}{}_p 
			+ \nabla_i \nabla^p T_{pj}{}^k
			+ \nabla^\qbar \hat{R}_{i \qbar j}{}^k.
\end{align*}
Hence $S$ satisfies the following evolution equation 
\begin{align}
	\left(\ddt - \Delta\right) S = {}& - |\nablabar \Psi|^2 - |\nabla \Psi|^2 \nonumber\\
		&+ \left( \nabla_r T_{\qbar}{}^{i \qbar} + \nabla_{\qbar} T^{\qbar}{}_r{}^i 
			\right) \Psi_{i j}{}^k \Psi^{r j}{}_k \nonumber
		+ \left( \nabla_r T_{\qbar}{}^{j \qbar} 
			+ \nabla_{\qbar} T^{\qbar}{}_{r}{}^j
			\right) \Psi_{ij}{}^k \Psi^{ir}{}_k	\nonumber \\
		& - \left( \nabla_k T_{\qbar}{}^{r \qbar} + \nabla_{\qbar} T^{\qbar}{}_k{}^r 
			\right) \Psi_{ij}{}^k \Psi^{ij}{}_r \nonumber\\
		&- 2 \mathrm{Re} \left[ \left(
			\nabla_i \nabla_j T^{p k}{}_p 
			+ \nabla_i \nabla_{\qbar} T^{\qbar}{}_j{}^k 
			- T_{ip}{}^r R_r{}^p{}_j{}^k 
			+ g^{\qbar p} \nabla_p \hat{R}_{i\qbar j}{}^k
			\right) \Psi^{ij}{}_k \right].
\label{evolution of S}
\end{align}
There are similar calculations to (\ref{evolution of S}) in the literature which generalize Calabi's argument \cite{Ca, Y}: in the elliptic Hermitian case \cite{Che, ZZ}; in  the case of the K\"ahler-Ricci flow (see also \cite{ShW}) in \cite{Cao, PSS}; and in other settings  \cite{TWY, To, ST}.

For the remainder of this section we will write $C$ for a constant of the form $CN^{\alpha}$ for $C$ and $\alpha$ depending only on $\hat{g}$.  Our goal is to show that $S \le C/r^2$.  The constant $C$ will be used repeatedly and may change from line to line, and we may at times use $C'$ or $C_1$ etc.

We would like to bound the right-hand side of (\ref{evolution of S}).  First, from (\ref{T = T hat}) and (\ref{nabla - nablahat = Psi}) we have, for example, 
\begin{equation}
	\nabla_\abar T_{ij}{}^k = g^{\lbar k} ( \hat \nabla_\abar \hat T_{ij\lbar} 
		- \Psi_{\abar\lbar}{}^\rbar \hat T_{ij\rbar} ).
\label{nabla T}
\end{equation}
This and similar calculations show that the second and third lines of (\ref{evolution of S}) can be bounded by $C(S^{3/2} + 1)$.  Next we address the terms in the last line of the evolution equation for $S$.  
\begin{itemize}

\item  
 Building on (\ref{nabla T}) we find
\begin{align} \nonumber
	\nabla_a \nabla_b T_{\ibar \jbar}{}^\kbar 
		={} & g^{\ov{k}l} \left( \nabla_a (\hat{\nabla}_b \hat{T}_{\ov{i} \ov{j}l} - \Psi_{bl}{}^r \hat{T}_{\ov{i}\ov{j} r} ) \right) \\
		 = {} & g^{ \kbar l} \left( \hat{\nabla}_a \hat{\nabla}_b \hat{T}_{\ibar\jbar l} 
		- \Psi_{a b}{}^r \hat{\nabla}_r \hat{T}_{\ibar \jbar l} 
		- \Psi_{a l}{}^r \hat{\nabla}_b \hat{T}_{\ibar \jbar r} \nonumber \right. \\
	&\left. - (\nabla_a \Psi_{bl}{}^r) \hat{T}_{\ov{i} \ov{j}r} 
		- \Psi_{bl}{}^r  \hat\nabla_a \hat{T}_{\ov{i} \ov{j} r} + \Psi_{b l}{}^{r} \Psi_{ar}{}^s \hat{T}_{\ov{i} \ov{j} s}  \right)
	\label{nabla nabla T},
\end{align}
and hence $|\nabla_i \nabla_j T^{p k}{}_p|$ can be bounded by $C (S+ | \nabla \Psi | + 1)$.

\item Similarly,
\begin{align} \nonumber
\nabla_a \nabla_{\ov{b}} T_{ij}{}^k ={} & g^{\ov{l}k} \nabla_a ( \hat{\nabla}_{\ov{b}} \hat{T}_{ij\ov{k}} - \Psi_{\ov{b} \ov{k}}{}^{\ov{q}} \hat{T}_{ij\ov{q}}) \\ \nonumber
={} & g^{\ov{l}k} \left( \hat{\nabla}_a \hat{\nabla}_{\ov{b}} \hat{T}_{ij\ov{k}} - \Psi_{ai}{}^p \hat{\nabla}_{\ov{b}} \hat{T}_{pj\ov{k}} - \Psi_{aj}{}^p \hat{\nabla}_{\ov{b}} \hat{T}_{ip\ov{k}} - (\nabla_a \Psi_{\ov{b} \ov{k}}{}^{\ov{q}}) \hat{T}_{ij\ov{q}} \right. \\
& \left. - \Psi_{\ov{b} \ov{k}}{}^{\ov{q}} ( \hat{\nabla}_a \hat{T}_{ij\ov{q}} - \Psi_{ai}{}^p \hat{T}_{pj \ov{q}} - \Psi_{aj}{}^p \hat{T}_{ip\ov{q}}) \right),
\end{align}
and so $|\nabla_i \nabla_{\qbar} T^{\qbar}{}_j{}^k|$ can be bounded by $C (S+ | \ov{\nabla} \Psi | + 1)$.

\item 
Next, using (\ref{T = T hat}) and (\ref{nablabar Psi}):  
\[
	T_{ip}{}^r R_r{}^p{}_j{}^k = g^{\sbar r}g^{\qbar p} \hat T_{ip\sbar} 
		\left( \hat R_{r\qbar j}{}^k - \nabla_\qbar \Psi_{rj}{}^k \right),
\]
so we can bound $|T_{ip}{}^r R_r{}^p{}_j{}^k|$ by $C(|\nablabar\Psi|+1)$.  

\item 
Finally, compute
\begin{align*}
	\nabla_p \hat R_{i\qbar j}{}^k 	& = \hat\nabla_p \hat R_{i\qbar j}{}^k - \Psi_{pi}{}^r \hat R_{r\qbar j}{}^k 
		- \Psi_{pj}{}^r \hat R_{i\qbar r}{}^k
		+ \Psi_{pr}{}^k \hat R_{i\qbar j}{}^r.
\end{align*}
So $|g^{\qbar p} \nabla_p \hat{R}_{i\qbar j}{}^k|$ can be bounded by $C(S^{1/2}+1)$. 
\end{itemize} 
Putting this all together  we arrive at the bound
\begin{equation}
	\left( \ddt - \Delta \right) S \le C (  S^{3/2} +  1 ) 
		- \frac12 ( |\nablabar\Psi|^2 + |\nabla\Psi|^2 ).
	\label{bound on evolution of S}
\end{equation}

We note here the bounds:
\begin{align}
	&\left| \nabla \tr \right|^2 \le C S 
		\label{bound grad trace g} \\
	&\left| \nabla S \right|^2 \le 2 S ( |\nablabar\Psi|^2+|\nabla\Psi|^2 ) .
		\label{bound grad S}
\end{align}
The first follows from 
$ \nabla_p \left( \hat{g}^{\jbar i} g_{i \jbar} \right)
	= \hat\nabla_p \left( \hat{g}^{\jbar i} g_{i \jbar} \right) 
	= \hat{g}^{\jbar i} \hat{\nabla}_p g_{i \jbar}  $
and the second follows from
$
	| \nabla S |^2 =  \bigl|\nabla |\Psi|^2\bigr| \, \bigl|\nablabar |\Psi|^2\bigr|
		\le 2 |\Psi|^2 \, ( |\nabla\Psi|^2 +    | \nablabar\Psi|^2)
$. 
Furthermore from \cite[Proposition 3.1]{TW1} (see also \cite{Che} in the elliptic case), we also have the following evolution equation for $\tr$:
\begin{align}
	\left( \ddt - \Delta \right) & \tr \ =\
		-g^{\jbar p} g^{\qbar i} \hat\nabla_k g_{i\jbar} \hat\nabla^k g_{p\qbar} 
		-2\operatorname{Re}\left(g^{\jbar i} \hat T_{ki}{}^p \hat\nabla^k g_{p\jbar} \right)
		\nonumber\\
	& 	+ g^{\jbar i} \left( \hat\nabla_i \hat T_{\jbar}{}^{k \qbar} 
		- \hat R_{i}{}^{k\qbar}{}_{\jbar} \right) g_{k\qbar} 
		- g^{\jbar i} \left( 
		\hat\nabla_i \hat T_{\jbar\qbar}{}^{\qbar} 
		+ \hat\nabla^k \hat T_{ik\jbar}
		\right) 
		\nonumber\\
	&  + g^{\jbar i} \hat T_{\jbar}{}^{k\qbar} \hat T_{ik}{}^p (\hat g - g)_{p\qbar} .
\label{evolution of trace g}
\end{align}
\newcommand{\mycomment}[1]{}%
\mycomment{
\begin{align}
	\left( \ddt - \Delta \right) & \tr \ =\
		-g^{\jbar p} g^{\qbar i} \hat\nabla_k g_{i\jbar} \hat\nabla^k g_{p\qbar} 
		\nonumber\\
	& -2\operatorname{Re}\left(g^{\jbar i} \hat T_{ki}{}^p \hat\nabla^k g_{p\jbar} \right)
		- g^{\jbar i} \hat T_{ik}{}^p \hat T_{\jbar}{}^{k \qbar} g_{p\qbar}
		+ g^{\jbar i} \left( \hat\nabla_i \hat T_{\jbar}{}^{k \qbar} 
		- \hat R_{i}{}^{k\qbar}{}_{\jbar} \right) g_{k\qbar} 
		\nonumber\\
	&  - g^{\jbar i} \left( 
		\hat\nabla_i \hat T_{\jbar\qbar}{}^{\qbar} 
		+ \hat\nabla^k \hat T_{ik\jbar}
		\right) 
		+ g^{\jbar i} \hat T_{\jbar}{}^{k\qbar} \hat T_{ik\qbar} .
\label{evolution of trace g}
\end{align}
}%
(Here $\hat\nabla^k = \hat g^{\lbar k}\hat\nabla_\lbar$ and we have raised indices on the tensor $\widehat\Rm$ using $\hat g$).
This generalizes the second order evolution inequality for the K\"ahler-Ricci flow \cite{Cao} (cf. \cite{Y,A}).
Hence we have the estimate 
\begin{equation}
	\left( \ddt - \Delta \right) \tr \le - \frac{S}{C_0} + C (S^{1/2}+1),
\label{bound on evolution of trace g}
\end{equation}
for a uniform positive constant $C_0$ (in fact we can take $C_0=N$).

We now would like to show that the evolution inequalities (\ref{bound on evolution of S}, \ref{bound on evolution of trace g}) imply a uniform bound on $S = |\hat \nabla g |^2$ on $\overline{B_{r/2}} \times [0,T]$.   Choose a smooth cutoff function $\rho$ which is supported in $B_r$ and is identically 1 on $\ov{B_{r/2}}$.  We may assume that 
$ |\nabla \rho|^2, |\Delta \rho|$ are bounded by $C /r^2$.  Let $K$ be a large uniform constant, to be specified later, which is at least large enough so that $$\frac K2 \le K - \tr \le K.$$  Let $A$ denote another  large positive constant to be specified later.  We will use a maximum principle argument with the function (cf. \cite{CY})
$$f = \rho^2 \frac{S}{K-\tr} + A \tr$$ to show that $S$ is bounded on $B_{r/2}$.  

Suppose that the maximum of $f$ on $\ov{B_r} \times [0,T]$ occurs at a point $(x_0, t_0)$.  We assume for the moment that  $t_0>0$ and that $x_0$ does not lie in the boundary of $\ov{B}_r$.
We wish to show that at $(x_0, t_0)$, $S$ is bounded from above by a uniform constant $C$.  Hence we 
 may assume without loss of generality that  $S>1$ at $(x_0, t_0)$.  In particular, we have
\begin{equation} \label{evolvestr}
	\left( \ddt - \Delta \right) S \le C   S^{3/2}  
		- \frac12 ( |\nablabar\Psi|^2 + |\nabla\Psi|^2 ),
		\quad 
		\left( \ddt - \Delta\right) \tr 
		\le -\frac{S}{2C_0} + 	C.
\end{equation}

We compute at $(x_0, t_0)$,
\begin{align*}
\left( \ddt - \Delta \right) f = {} & A \left( \ddt - \Delta \right) \tr + (-\Delta (\rho^2)) \frac{S}{K-\tr} + \rho^2 \frac{S}{(K-\tr)^2} \left( \ddt - \Delta \right) \tr  \\
& + \rho^2 \frac{1}{K-\tr} \left( \ddt - \Delta \right) S - 4 \textrm{Re} \left[ \rho \frac{S}{(K-\tr)^2} \nabla \tr \cdot \ov{\nabla} \rho \right] \\
& - 4 \textrm{Re} \left[ \rho \frac{1}{K-\tr} \nabla \rho \cdot \ov{\nabla} S \right] - 2 \textrm{Re} \left[ \rho^2 \frac{1}{(K-\tr)^2} \nabla \tr \cdot \ov{\nabla} S \right]  \\
& - \frac{2\rho^2 S}{(K-\tr)^3} | \nabla \tr|^2.
\end{align*}
But since a maximum occurs at $(x_0, t_0)$ we have $\ov{\nabla} f=0$ at this point, and hence
$$2\rho \ov{\nabla}\rho \frac{S}{K-\tr} + \rho^2 \frac{\ov{\nabla} S}{K-\tr} + \rho^2 \frac{S \ov{\nabla} \tr}{(K-\tr)^2} + A \ov{\nabla} \tr =0.$$
Then at $(x_0, t_0)$,
\begin{align*}
\left( \ddt - \Delta \right) f = {} & A \left( \ddt - \Delta \right) \tr + (-\Delta (\rho^2)) \frac{S}{K-\tr} + \rho^2 \frac{S}{(K-\tr)^2} \left( \ddt - \Delta \right) \tr \\
& + \rho^2 \frac{1}{K-\tr} \left( \ddt - \Delta \right) S - 4 \textrm{Re} \left[ \rho  \frac{1}{K-\tr} \nabla \rho \cdot \ov{\nabla} S \right] + \frac{2A | \nabla \tr|^2}{K-\tr}.
\end{align*}
Making use of (\ref{bound grad trace g}, \ref{bound grad S}, \ref{evolvestr}) and Young's inequality, we obtain at $(x_0, t_0)$,
\begin{align*}
0 \le \left( \ddt - \Delta \right) f \le & \left( -\frac{A}{2C_0} S + CA \right) + \left( \frac{CS}{r^2K} \right) + \left( - \frac{\rho^2}{2K^2 C_0} S^2 + \frac{C\rho^2}{K^2} S \right) \\
& + \left( -\frac{\rho^2}{2K} ( | \ov{\nabla} \Psi|^2 + | \nabla \Psi|^2 ) + \frac{\rho^2}{4K^2C_0} S^2 + C \rho^2 S\right) \\
& + \left( \frac{\rho^2}{4K} ( | \ov{\nabla} \Psi|^2 + | \nabla \Psi|^2 ) + \frac{C}{Kr^2} S \right) + \frac{CA}{K} S \\
& \le - \frac{A}{2C_0} S + CA + \frac{C'}{r^2} S + \frac{CA}{K} S.
\end{align*}
Now pick $K\ge4C_0 C$ so that  at $(x_0, t_0)$,
$$0 \le - \frac{A}{4C_0}S + CA + \frac{C'}{r^2} S.$$
Then choose $A= \frac{8C'C_0}{r^2}$ so that at $(x_0, t_0)$,
$$\frac{C'}{r^2} S \le CA,$$
giving a uniform upper bound for $S$.  It follows that $f$ is bounded from above by $Cr^{-2}$ for a uniform $C$.  Hence $S$ on $\ov{B_{r/2}}$ is bounded above by $Cr^{-2}$.

It remains to deal with the cases when $t_0=0$ or $x_0$ lies on the boundary of $\ov{B_r}$.  In either case we have $f(x_0, t_0) \le  A \tr (x_0, t_0) \le Cr^{-2}$ and the same bound holds.

\begin{remark} \label{remark1} \emph{
Tracing through the argument, one can see that the constants only depend on uniform bounds for the torsion and curvature of $\hat{g}$, and one and two derivatives (with respect to $\hat{\nabla}$ or $\ov{\hat{\nabla}}$) of torsion and one derivative of curvature.}
%
\end{remark}

\section{Local curvature bound}

In this section we prove part (ii) of Theorem \ref{maintheorem}.
As in the previous section, we write $C$ for a constant of the form $CN^{\gamma}$ for some uniform $C, \gamma$. We compute in the ball $\ov{B_{r/2}}$ on which we already have the bound $S\le C/r^2$.

Let $\Delta_{\mathbb R} = \frac12 g^{\qbar p}( \nabla_p \nabla_\qbar + \nabla_\qbar \nabla_p )$.  First we need an evolution equation for the curvature tensor.  
We begin with 
\[
	\ddt R_{i \jbar k}{}^l = \ddt \left( - \partial_\jbar \Gamma_{ik}^l \right) 
	 = - \partial_\jbar \ddt \left( \Gamma_{ik}^l \right) 
	 = - \partial_\jbar ( - \nabla_i (R^C)_k{}^l ) 
	 = \nabla_\jbar \nabla_i R_k{}^l{}_p{}^p 
\]
and therefore, 
\begin{equation}
	\ddt R_{i \jbar k \lbar} = - R_{q \lbar p}{}^p R_{i \jbar k}{}^q 
		+ \nabla_\jbar \nabla_i R_{k \lbar p}{}^p.
\end{equation}
Now, computing in coordinates where $g$ is the identity, we find
\begin{align*}
	\Delta_{\mathbb R} R_{i \jbar k \lbar} 	=& \frac12 ( \nabla_p \nabla_\pbar + \nabla_\pbar \nabla_p ) R_{i \jbar k \lbar} \\
	=& \nabla_p \nabla_\pbar R_{i \jbar k \lbar} + \frac12 (
		R_{p \pbar i \qbar} R_{q \jbar k \lbar} 
		- R_{p \pbar q \jbar} R_{i \qbar k \lbar}
		+ R_{p \pbar k \qbar} R_{i \jbar q \lbar}
		- R_{p \pbar q \lbar} R_{i \jbar k \qbar}
		) \\
	=& \nabla_p ( \nabla_\jbar R_{i \pbar k \lbar} - T_{\pbar \jbar q} R_{i \qbar k \lbar} ) 
		+ \frac12 (
		R_{p \pbar i \qbar} R_{q \jbar k \lbar} 
		- R_{p \pbar q \jbar} R_{i \qbar k \lbar}
		+ R_{p \pbar k \qbar} R_{i \jbar q \lbar}
		- R_{p \pbar q \lbar} R_{i \jbar k \qbar}
		) \\
	=& \nabla_\jbar \nabla_p R_{i \pbar k \lbar} 
		- R_{p \jbar i \qbar} R_{q \pbar k \lbar}
		+ R_{p \jbar q \pbar} R_{i \qbar k \lbar}
		- R_{p \jbar k \qbar} R_{i \pbar q \lbar}
		+ R_{p \jbar q \lbar} R_{i \pbar k \qbar} 
		- \nabla_p ( T_{\pbar \jbar q} R_{i \qbar k \lbar} )  \\
	& + \frac12 (
		R_{p \pbar i \qbar} R_{q \jbar k \lbar} 
		- R_{p \pbar q \jbar} R_{i \qbar k \lbar}
		+ R_{p \pbar k \qbar} R_{i \jbar q \lbar}
		- R_{p \pbar q \lbar} R_{i \jbar k \qbar}
		) \\
	=& \nabla_\jbar ( \nabla_i R_{p \pbar k \lbar} - T_{p i \qbar} R_{q \pbar k \lbar} ) 
		- R_{p \jbar i \qbar} R_{q \pbar k \lbar}
		+ R_{p \jbar q \pbar} R_{i \qbar k \lbar}
		- R_{p \jbar k \qbar} R_{i \pbar q \lbar}
		+ R_{p \jbar q \lbar} R_{i \pbar k \qbar} \\
	& - \nabla_p ( T_{\pbar \jbar q} R_{i \qbar k \lbar} )
		+ \frac12 (
		R_{p \pbar i \qbar} R_{q \jbar k \lbar} 
		- R_{p \pbar q \jbar} R_{i \qbar k \lbar}
		+ R_{p \pbar k \qbar} R_{i \jbar q \lbar}
		- R_{p \pbar q \lbar} R_{i \jbar k \qbar}
		) \\
	=& \nabla_\jbar \nabla_i ( R_{k \lbar p \pbar} 
		- \nabla_p T_{\pbar \lbar k} 
		- \nabla_\lbar T_{p k \pbar} 
		)
		- \nabla_\jbar ( T_{p i \qbar} R_{q \pbar k \lbar} ) \\
	& - R_{p \jbar i \qbar} R_{q \pbar k \lbar}
		+ R_{p \jbar q \pbar} R_{i \qbar k \lbar}
		- R_{p \jbar k \qbar} R_{i \pbar q \lbar}
		+ R_{p \jbar q \lbar} R_{i \pbar k \qbar} \\
	& - \nabla_p ( T_{\pbar \jbar q} R_{i \qbar k \lbar} ) 
		+ \frac12 (
		R_{p \pbar i \qbar} R_{q \jbar k \lbar} 
		- R_{p \pbar q \jbar} R_{i \qbar k \lbar}
		+ R_{p \pbar k \qbar} R_{i \jbar q \lbar}
		- R_{p \pbar q \lbar} R_{i \jbar k \qbar}
		) . 
\end{align*}
Hence
\begin{align} \nonumber
	\left( \ddt - \Delta_{\mathbb R} \right) R_{i \jbar k \lbar} 
		=& - R_{q \lbar p \pbar} R_{i \jbar k \qbar}  
		+ R_{p \jbar i \qbar} R_{q \pbar k \lbar}
		- R_{p \jbar q \pbar} R_{i \qbar k \lbar}
		+ R_{p \jbar k \qbar} R_{i \pbar q \lbar}
		- R_{p \jbar q \lbar} R_{i \pbar k \qbar} \\ \nonumber
	& - \frac12 (
		R_{p \pbar i \qbar} R_{q \jbar k \lbar} 
		- R_{p \pbar q \jbar} R_{i \qbar k \lbar}
		+ R_{p \pbar k \qbar} R_{i \jbar q \lbar}
		- R_{p \pbar q \lbar} R_{i \jbar k \qbar}
		) \\ \label{evRm}
	& + \nabla_p ( \hat{T}_{\pbar \jbar q} R_{i \qbar k \lbar} ) 
		+ \nabla_\jbar ( \hat{T}_{p i \qbar} R_{q \pbar k \lbar} )
		+ \nabla_\jbar \nabla_i ( \nabla_p T_{\pbar \lbar k} 
			+ \nabla_\lbar T_{p k \pbar} ).
\end{align}
To estimate this, we first compute
$$\nabla_p (\hat{T}_{\ov{p} \ov{j} q} R_{i\ov{q}k \ov{l}}) = (\hat{\nabla}_p \hat{T}_{\ov{p} \ov{j} q} - \Psi_{pq \ov{r}} \hat{T}_{\ov{p}\ov{j}r}) R_{i\ov{q}k\ov{l}} + \hat{T}_{\ov{p}\ov{j} q} \nabla_p R_{i\ov{q}k\ov{l}},$$
and this is bounded by $C(|\textrm{Rm}|/r + | \nabla \textrm{Rm}|)$.   Using the fact that $R_{i\ov{j}k}{}^l = - \nabla_{\ov{j}} \Psi_{ik}{}^l + \hat{R}_{i\ov{j}k}{}^l$ we have
\begin{equation} \label{bdRmPsi}
|\textrm{Rm}| \le |\ov{\nabla} \Psi|+C,
\end{equation}
and hence
\begin{equation} \label{bdThatR}
|\nabla_p (\hat{T}_{\ov{p} \ov{j} q} R_{i\ov{q}k \ov{l}}) | \le C \left(  | \nabla \textrm{Rm}| + \frac{|\ov{\nabla}\Psi|}{r} + \frac{1}{r} \right).
\end{equation}
Similarly for the term $ \nabla_\jbar ( \hat{T}_{p i \qbar} R_{q \pbar k \lbar} )$. 

The last  two terms of (\ref{evRm}) involve three derivatives of torsion.  We claim that
\begin{equation} \label{claimT}
	| \nablabar \nabla \nabla \overline T |, \ | \nablabar \nabla \nablabar T |
		\le C \left(| \nabla \textrm{Rm}| + \frac{|\nabla\Psi|+|\nablabar\Psi|}{r} + \frac{1}{r^3} \right).
\end{equation}
Indeed, applying $\nabla_{\ov{c}}$ to (\ref{nabla nabla T}), we have
\begin{align} \nonumber
\nabla_{\ov{c}} \nabla_a \nabla_b T_{\ov{i} \ov{j}}{}^{\ov{k}} = {} & g^{\ov{k}l} \left( \hat{\nabla}_{\ov{c}} \hat{\nabla}_a \hat{\nabla}_b \hat{T}_{\ov{i}\ov{j}l} - \Psi_{\ov{c}\ov{i}}{}^{\ov{q}} \hat{\nabla}_a \hat{\nabla}_b \hat{T}_{\ov{q} \ov{j} l} - \Psi_{\ov{c} \ov{j}}{}^{\ov{q}} \hat{\nabla}_a \hat{\nabla}_b \hat{T}_{\ov{i}\ov{q}l} \right. \\ \nonumber
& - \nabla_{\ov{c}} (\Psi_{ab}{}^r \hat{\nabla}_r \hat{T}_{\ov{i} \ov{j}l}) - \nabla_{\ov{c}} (\Psi_{al}{}^r \hat{\nabla}_b \hat{T}_{\ov{i}\ov{j}r}) - \nabla_{\ov{c}} (\Psi_{bl}{}^r \hat{\nabla}_a \hat{T}_{\ov{i}\ov{j}r}) \\ \label{3derivsT}
&\left. - \nabla_{\ov{c}} (\nabla_a \Psi_{bl}{}^r \hat{T}_{\ov{i} \ov{j}r}) + \nabla_{\ov{c}} (\Psi_{bl}{}^r \Psi_{ar}{}^s \hat{T}_{\ov{i} \ov{j}s}) \right).
\end{align}
The first three terms on the right hand side of (\ref{3derivsT}) are bounded by $C (\sqrt{S}+1)$ and hence by $C/r$.  Next compute
\begin{align} \label{1}
\nabla_{\ov{c}} (\Psi_{ab}{}^r \hat{\nabla}_r \hat{T}_{\ov{i}\ov{j}l}) = (\nabla_{\ov{c}} \Psi_{ab}{}^r) \hat{\nabla}_r \hat{T}_{\ov{i}\ov{j}l} + \Psi_{ab}{}^r \hat{\nabla}_{\ov{c}} \hat{\nabla}_r \hat{T}_{\ov{i} \ov{j}l} - \Psi_{ab}{}^r \Psi_{\ov{c}\ov{i}}{}^{\ov{q}} \hat{\nabla}_r \hat{T}_{\ov{q}\ov{j}l} - \Psi_{ab}{}^r \Psi_{\ov{c}\ov{j}}{}^{\ov{q}} \hat{\nabla}_r \hat{T}_{\ov{i}\ov{q}l},
\end{align}
which is bounded by $C | \ov{\nabla} \Psi| + C\sqrt{S}  + C S$ and hence by  $C(| \ov{\nabla} \Psi| + 1/r^2)$.  The same bound holds for the other two terms on the second line of (\ref{3derivsT}).

For the third line, compute
\begin{align*}
\nabla_{\ov{c}} (\nabla_a \Psi_{bl}{}^r \hat{T}_{\ov{i}\ov{j}r}) = {}& \left( \nabla_a \nabla_{\ov{c}} \Psi_{bl}{}^r + R_{a\ov{c}b}{}^p \Psi_{pl}{}^r + R_{a\ov{c}l}{}^p \Psi_{bp}{}^r - R_{a\ov{c}p}{}^r \Psi_{bl}{}^p\right) \hat{T}_{\ov{i}\ov{j}r} \\
& + (\nabla_a \Psi_{bl}{}^r) \left( \hat{\nabla}_{\ov{c}} \hat{T}_{\ov{i} \ov{j} r} - \Psi_{\ov{c} \ov{i}}{}^{\ov{q}} \hat{T}_{\ov{q} \ov{j}r} - \Psi_{\ov{c} \ov{j}}{}^{\ov{q}} \hat{T}_{\ov{i} \ov{q}r} \right),
\end{align*}
and using the fact that $\nabla_{\ov{c}} \Psi_{bl}{}^r = - R_{b\ov{c} l}{}^r + \hat{R}_{b\ov{c} l}{}^r$ we obtain
\begin{align*}
\nabla_{\ov{c}} (\nabla_a \Psi_{bl}{}^r \hat{T}_{\ov{i}\ov{j}r}) = {}& \big( - \nabla_a R_{b\ov{c} l}{}^r + \hat{\nabla}_a \hat{R}_{b\ov{c}l}{}^r - \Psi_{ab}{}^p \hat{R}_{p\ov{c}l}{}^r - \Psi_{al}{}^p \hat{R}_{b\ov{c}p}{}^r + \Psi_{ap}{}^r \hat{R}_{b\ov{c}l}{}^p  \\ 
& + R_{a\ov{c}b}{}^p \Psi_{pl}{}^r + R_{a\ov{c}l}{}^p \Psi_{bp}{}^r - R_{a\ov{c}p}{}^r \Psi_{bl}{}^p\big) \hat{T}_{\ov{i}\ov{j}r} \\
& + (\nabla_a \Psi_{bl}{}^r) \left( \hat{\nabla}_{\ov{c}} \hat{T}_{\ov{i} \ov{j} r} - \Psi_{\ov{c} \ov{i}}{}^{\ov{q}} \hat{T}_{\ov{q} \ov{j}r} - \Psi_{\ov{c} \ov{j}}{}^{\ov{q}} \hat{T}_{\ov{i} \ov{q}r} \right).
\end{align*}
It follows that 
\begin{equation} \label{2}
| \nabla_{\ov{c}} (\nabla_a \Psi_{bl}{}^r \hat{T}_{\ov{i}\ov{j}r})| \le C \left( | \nabla \textrm{Rm}| + \frac{| \textrm{Rm}|}{r} + \frac{|\nabla \Psi|}{r} +\frac{1}{r} \right).
\end{equation}
Finally,
\begin{align*}
\nabla_{\ov{c}} (\Psi_{bl}{}^r \Psi_{ar}{}^s \hat{T}_{\ov{i}\ov{j}s}) = {} & ( - R_{b\ov{c} l}{}^r + \hat{R}_{b\ov{c}l}{}^r) \Psi_{ar}{}^s \hat{T}_{\ov{i}\ov{j}s} + \Psi_{bl}{}^r ( -R_{a\ov{c} r}{}^s + \hat{R}_{a\ov{c}r}{}^s) \hat{T}_{\ov{i}\ov{j}s} \\
& + \Psi_{bl}{}^r \Psi_{ar}{}^s ( \hat{\nabla}_{\ov{c}} \hat{T}_{\ov{i} \ov{j}s} - \Psi_{\ov{c}\ov{i}}{}^{\ov{q}} \hat{T}_{\ov{q}\ov{j}s} - \Psi_{\ov{c}\ov{j}}{}^{\ov{q}} \hat{T}_{\ov{i}\ov{q}s}),
\end{align*}
giving
\begin{equation} \label{3}
|\nabla_{\ov{c}} (\Psi_{bl}{}^r \Psi_{ar}{}^s \hat{T}_{\ov{i}\ov{j}s})|  \le C\left( \frac{|\textrm{Rm}|}{r} + \frac{1}{r^3} \right).
\end{equation}
Putting together (\ref{3derivsT}, \ref{1}, \ref{2}, \ref{3}), and making use of (\ref{bdRmPsi}), we obtain
$$| \nablabar \nabla \nabla \overline T |
		\le C \left(| \nabla \textrm{Rm}| + \frac{|\nabla\Psi|+|\nablabar\Psi|}{r} + \frac{1}{r^3} \right),$$
and the bound for $| \nablabar \nabla \nablabar T |$ follows similarly.  This completes the proof of the claim (\ref{claimT}).

From (\ref{bdThatR}) and the claim we just proved,  since the first two lines of (\ref{evRm}) are of the order $|\textrm{Rm}|^2$, we have the bound
\begin{equation}
	\left| \left( \ddt - \Delta_{\mathbb R} \right) \Rm \right| 
		\le C \left( |\textrm{Rm}|^2 + | \nabla \textrm{Rm} |+  \frac{| \nabla \Psi|+ | \ov{\nabla} \Psi|}{r} + \frac{1}{r^3} \right).
\label{bound on evolution of curvature tensor}
\end{equation} 
	
Now
\begin{align}
	\left( \ddt - \Delta \right) | \Rm |^2 
		={} & g^{\jbar b} g^{\cbar k} g^{\lbar d} (R^C)^{\abar i} 
		R_{i \jbar k \lbar} \overline{R_{a \bbar c \dbar}} \nonumber\\ 
	&+ g^{\abar i} g^{\cbar k} g^{\lbar d} (R^C)^{\jbar b} 
		R_{i \jbar k \lbar} \overline{R_{a \bbar c \dbar}} \nonumber\\
	&+ g^{\abar i} g^{\jbar b} g^{\lbar d} (R^C)^{\cbar k} 
		R_{i \jbar k \lbar} \overline{R_{a \bbar c \dbar}} \nonumber\\
	&+ g^{\abar i} g^{\jbar b} g^{\cbar k} (R^C)^{\lbar d} 
		R_{i \jbar k \lbar} \overline{R_{a \bbar c \dbar}} \nonumber\\
	&+ 2 \mathrm{Re} \left[ 
		g^{\abar i} g^{\jbar b} g^{\cbar k} g^{\lbar d} 
		\left( (\ddt-\Delta_{\mathbb R}) R_{i \jbar k \lbar} \right)
		\overline{R_{a \bbar c \dbar}}
		\right] \nonumber\\
	&- 2 |\nabla \Rm|^2.
\label{evolution of Rm squared}
\end{align}
This together with (\ref{bound on evolution of curvature tensor}) and (\ref{bdRmPsi}) implies
\begin{align} \nonumber
	\left( \ddt - \Delta \right) | \Rm |^2 \le {} & C \left( |\textrm{Rm}|^2 +  |\textrm{Rm}|^3 + |\nabla \textrm{Rm}| \cdot |\textrm{Rm}| + \frac{(| \nabla \Psi|+ |\ov{\nabla}\Psi|) | \textrm{Rm}|}{r} + \frac{| \textrm{Rm}|}{r^3} \right) \\
		\le {} & C \left( |\Rm|^3+\frac{1}{r}  +  \frac{ |\nabla\Psi|^2+|\nablabar\Psi|^2}{r} + \frac{|\textrm{Rm}|}{r^3} \right) - |\nabla\Rm|^2.
\label{bound on evolution of Rm squared}
\end{align}

To show $|\Rm|^2$ is locally uniformly bounded we will use an argument similar to the previous section.  Let $\rho$ now denote a cutoff function which is identically 1 on $\ov{B_{r/4}}$, and supported in $B_{r/2}$.  From the previous section we know that $S$ is bounded by $C/r^2$ on $B_{r/2}$.  As before we can assume $|\nabla\rho|^2$ and $|\Delta\rho|$ are bounded by $C/r^2$.  Let $K= C_1/r^2$ where $C_1$ is a constant to be determined later, and is at least   large enough so that $\frac{K}2 \le K-S \le K$.  Let $A$ denote a constant to be specified later.  We will apply the maximum principle argument to the quantity 
$$f = \rho^2 \frac{|\Rm|^2}{K-S} + A S.$$  
As in the previous section, we calculate at a point $(x_0, t_0)$ where a maximum of $f$ is achieved, and we first assume that $t_0>0$ and that $x_0$ does not occur at the boundary of $\ov{B_{r/2}}$.  We use the fact that $\nabla f=0$ at this point, giving us
\begin{align*}
	\left( \ddt - \Delta \right)  f
		=& A ( \ddt - \Delta ) S 
	+ (- \Delta (\rho^2) ) \frac{ |\Rm|^2 }{K - S} + \rho^2 \frac{|\Rm|^2}{(K - S)^2} ( \ddt - \Delta ) S \\
	& + \rho^2 \frac1{K - S} ( \ddt - \Delta ) |\Rm|^2 
	- 4 \mathrm{Re} \left(
		\frac1{K - S} \rho \nabla \rho \cdot \nablabar |\Rm|^2 
		\right) + \frac{2A | \nabla S|^2}{K - S}.
		\end{align*}
Our goal is to show that at $(x_0, t_0)$, we have $|\textrm{Rm}|^2 \le C/r^4$.   Hence without loss of generality, we may assume that $1/r + | \textrm{Rm}|/r^3 \le C |\textrm{Rm}|^3$ and hence (\ref{bound on evolution of Rm squared}) becomes
$$\left( \ddt - \Delta \right) | \Rm |^2 
		\le C \left( |\Rm|^3  +  \frac{Q}{r} \right) - |\nabla\Rm|^2,
$$
where for convenience we are writing $Q = |\nabla\Psi|^2+|\nablabar\Psi|^2$.  For later purposes, recall  from (\ref{bdRmPsi}) that $|\textrm{Rm}|^2 \le Q+C$ and from (\ref{bound grad S}) that $|\nabla S|^2 \le 2SQ$.

 Also note that $\bigl|\nabla|\Rm|^2\bigr| \le 2 |\Rm| |\nabla\Rm|$.  
By (\ref{bound on evolution of S}) we find that on $B_{r/2}$ we have $$(\ddt - \Delta)S \le \frac{C}{r^3} - \frac12 Q.$$  Using these, we find at $(x_0, t_0)$,
\begin{align*}
\left( \ddt -\Delta \right) f \le {} & \left( \frac{CA}{r^3} - \frac{AQ}{2} \right) + \left( \frac{C|\textrm{Rm}|^2}{Kr^2} \right) + \left( \frac{C \rho^2 |\textrm{Rm}|^2}{K^2r^3}  - \frac{\rho^2 | \textrm{Rm}|^2 Q}{2K^2} \right)  \\
& + \left( \frac{C\rho^2 |\textrm{Rm}|^3}{K} + \frac{C\rho^2Q}{Kr} - \frac{\rho^2}{K} | \nabla \textrm{Rm} |^2 \right) + \left( \frac{\rho^2 | \nabla \textrm{Rm}|^2}{2K} + C \frac{| \textrm{Rm}|^2}{Kr^2} \right) + \left( \frac{8ASQ}{K} \right).
\end{align*}
First choose $C_1$ in the definition of $K$ to be sufficiently large so that 
$$\frac{8ASQ}{K} \le \frac{AQ}{4},$$
where we use the fact that $S\le C/r^2$.  Next observe that
$$\frac{C \rho^2 |\textrm{Rm}|^3}{K} \le \frac{\rho^2 |\textrm{Rm}|^2 Q}{2K^2} + C' \rho^2 |\textrm{Rm}|^2,$$
and hence
\begin{align*}
\left( \ddt - \Delta \right) f \le & \frac{CA}{r^3} - \frac{AQ}{4} + C''Q + C.
\end{align*}
Now we may choose $A$ sufficiently large so that $A \ge 8C''$ and we obtain at $(x_0, t_0)$,
$$Q \le \frac{C}{r^3},$$
which implies that $|\textrm{Rm}|^2 \le C/r^3$ at this point.  It follows that at $(x_0, t_0)$, $f$ is bounded from above by $C/r^2$.  The same bound holds if $x_0$ lies in the boundary of $\ov{B_{r/2}}$ or if $t_0=0$.  Hence on $\ov{B_{r/4}}$ we obtain
$$|\textrm{Rm}|^2 \le \frac{C}{r^4},$$
as required.  This completes the proof of Theorem \ref{maintheorem}.

\begin{remark} \label{remark2} \emph{In addition to the dependence discussed in Remark \ref{remark1}, the constants also depend on \emph{three} derivatives of the torsion of $\hat{g}$, with respect to $\hat{\nabla}$ or $\ov{\hat{\nabla}}$}.

\end{remark}

\section{Higher order estimates}

In this last section, we prove Corollary \ref{corollary} by establishing the estimates for $|\hat{\nabla}_{\mathbb{R}}^m g|^2_{\hat{g}}$ for $m=2, 3, \ldots$.  For this part, we essentially follow the method of Gill \cite{G} (cf. \cite{Chau, CK, CLN, PSSW} in the K\"ahler case), but since the setting here is slightly more general, we briefly outline the argument.  In this section, we say that a quantity is uniformly bounded if it can be bounded by $CN^{\alpha} r^{-\gamma}$ for uniform $C, \alpha, \gamma$.

We work on the ball $B_{r/4}$, and assume the bounds established in Theorem \ref{maintheorem}.
As in \cite{TW1}, define reference tensors $(\hat{g}_t)_{i\ov{j}} = \hat{g}_{i\ov{j}} - t \hat{R}^C_{i\ov{j}}$, where $\hat{R}^C_{i\ov{j}}$ is the Chern-Ricci curvature of $\hat{g}$.  
For each fixed $x\in M$, let $\varphi=\varphi(x,t)$ solve 
$$\frac{\partial \varphi}{\partial t} = \log \frac{\det g(t)}{\det \hat{g}}, \quad \varphi|_{t=0} =0.$$
Then $g_{i\ov{j}} = (\hat{g}_t)_{i\ov{j}} +\partial_i \partial_{\ov{j}} \varphi$ is the solution of the Chern-Ricci flow starting at $\hat{g}$.

Consider the  first order differential operator $D = \frac{\partial}{\partial x^{\gamma}}$, where $x^{\gamma}$ is a real coordinate.
Applying $D$ to the equation for $\varphi$, we have
\begin{equation*} 
\frac{\partial}{\partial t} (D\varphi) = g^{\ov{j}i} D g_{i\ov{j}} - \hat{g}^{\ov{j}i} D (\hat{g})_{i \ov{j}}
 = g^{\ov{j}i} \partial_i \partial_{\ov{j}} (D\varphi) + g^{\ov{j}i} D (\hat{g}_t)_{i\ov{j}} - \hat{g}^{\ov{j}i} D \hat{g}_{i \ov{j}}.
 \end{equation*}
Hence, working in real coordinates, the function $u=D(\varphi)$ satisfies a linear parabolic PDE of the form
\begin{equation} \label{linearPDE}
\partial_t u = a^{\alpha \beta} \partial_{x^{\alpha}} \partial_{x^{\beta}} u + f,
\end{equation}
where $A=(a^{\alpha \beta})$ is a real $2n \times 2n$ positive definite symmetric matrix whose largest and smallest eigenvalues $\Lambda$ and $\lambda$ satisfy 
\begin{equation} \label{evalues}
C^{-1} \le \lambda \le \Lambda \le C,
\end{equation}
for a uniform positive constant $C$.

Moreover, the entries of $A$ are uniformly bounded in the $C^{\delta/2,\delta}$ parabolic norm for $0< \delta<1$.  Indeed our Calabi-type estimate from part (i) of Theorem \ref{maintheorem},
$$| \hat{\nabla}g|^2 \le C,$$
implies that the Riemannian metric  $g_R$ associated to $g$ is bounded in the $C^1$ norm in the space direction.    On the other hand, 
$$\ddt{} g_{i\ov{j}} = - R^{C}_{i\ov{j}} = - g^{\ov{l}k} R_{i\ov{j}k\ov{l}}.$$
From the curvature bound of Theorem \ref{maintheorem}, we know that $ g^{\ov{l}k} R_{i\ov{j}k\ov{l}}$ is uniformly bounded for any fixed $i,j$.
It follows that $\ddt{} (g_R)_{\alpha \beta}$ is also uniformly bounded for any fixed $\alpha, \beta$.   Thus we see that each entry $a^{\alpha \beta}$ in the matrix $A$ has uniform bounds in one space and one time derivative.  This implies that $a^{\alpha \beta}$ is uniformly bounded in the $C^{\delta/2, \delta}$ parabolic norm for any $0< \delta<1$.  

Next, note that $u= \frac{\partial \varphi}{\partial x^{\gamma}}$ in (\ref{linearPDE}) is bounded in the $C^0$ norm since $g(t)$ is uniformly bounded and hence $| \ddbar \varphi|_{C^0}$ is uniformly bounded.   Moreover, $f$ in (\ref{linearPDE}) is uniformly bounded in the $C^{\delta/2, \delta}$ norm.

We can then apply Theorem 8.11.1 in \cite{Kr} to (\ref{linearPDE}) to see that $u$ is bounded in the parabolic $C^{1+\delta/2, 2+\delta}$ norm on a slightly smaller parabolic domain: $[\ve', T] \times B_{r'}$ for any $\ve'$ and $r'$ with $0 < \ve' < \ve$ and $r/8 < r' <r/4$.
Tracing through the argument in \cite{Kr}, one can check that the estimates we obtain indeed are of the desired form.

Now apply $D$ to the equality $g_{i\ov{j}}(t) = (\hat{g}_t)_{i\ov{j}} + \partial_i \partial_{\ov{j}} \varphi$.  We get
$$Dg_{i\ov{j}} = D(\hat{g}_t)_{i\ov{j}}  + \partial_i \partial_{\ov{j}} u,$$
where we recall that $D= \partial/\partial x^{\gamma}$ for some $\gamma$.  Since we have bounds for $u$ in $C^{1+\delta/2, 2+\delta}$ this implies that $\partial_i \partial_{\ov{j}} u$ is bounded in $C^{\delta/2, \delta}$.  Since $ D(\hat{g}_t)_{i\ov{j}} $ is uniformly bounded in all norms  we get that $D g_{i\ov{j}}$ is uniformly bounded in $C^{\delta/2, \delta}$ for all $i,j$.  Since $D= \partial/\partial x^{\gamma}$ and $\gamma$ was an arbitrary index, it follows that $\partial_{\gamma} a^{\alpha \beta}$ is uniformly bounded in $C^{\delta/2, \delta}$ for all $\alpha, \beta, \gamma$.  We have  a similar estimate for $\partial_{\gamma} f$. Now apply Theorem 8.12.1 in \cite{Kr} (with $k=1$) to see that, for any $\alpha$,  $\partial_{\alpha} u$ is uniformly bounded in $C^{1+\delta/2, 2+\delta}$ on a slightly smaller parabolic domain.  This means that $D^{{\alpha}} \varphi$ is uniformly bounded in $C^{1+\delta/2, 2+\delta}$ for any multi-index $\alpha\in \mathbb{R}^{2n}$ with $| \alpha| \le 2$.  

We can then iterate this procedure and obtain the required $C^k$ bounds for $g(t)$ for all $k$. This completes the proof of the corollary.

\begin{remark} \emph{In \cite{ShW}, we showed how to obtain higher derivative estimates for curvature using simple maximum principle arguments (following \cite{H, Shi}).  However, in the case of the Chern-Ricci flow, there are difficulties in using this approach because of torsion terms that need to be controlled.  An alternative method to proving the estimates in this section may be to generalize the work of Gill on the K\"ahler-Ricci flow \cite{G2}.  This could give an ``elementary'' maximum principle proof, but the technical difficulties in carrying this out seem to be substantial.}
\end{remark}


\end{document}